\documentclass[12pt]{article}%
\usepackage{graphicx}
\usepackage[intlimits]{amsmath}
\usepackage{amsfonts}
\usepackage{amssymb}%
 \makeatletter 
\newtheorem{theorem}{Theorem}

\newenvironment{proof}[1][Proof]{\noindent{\textbf {#1}  }}  {\hfill$\Box$\bigskip}

\begin{document}

\title{Graphs and matrices with maximal energy}
\author{Vladimir Nikiforov\\Department of Mathematical Sciences, University of Memphis, \\Memphis TN 38152, USA, e-mail:\textit{ vnikifrv@memphis.edu}}
\maketitle

\begin{abstract}
Given a complex $m\times n$ matrix $A,$ we index its singular values as
$\sigma_{1}\left(  A\right)  \geq\sigma_{2}\left(  A\right)  \geq...$ and call
the value $\mathcal{E}\left(  A\right)  =\sigma_{1}\left(  A\right)
+\sigma_{2}\left(  A\right)  +...$ the \emph{energy }of $A,$ thereby extending
the concept of graph energy, introduced by Gutman. Koolen and Moulton proved
that $\mathcal{E}\left(  G\right)  \leq\left(  n/2\right)  \left(  1+\sqrt
{n}\right)  $ for any graph $G$ of order $n$ and exhibited an infinite family
of graphs with $\mathcal{E}\left(  G\right)  =\left(  v\left(  G\right)
/2\right)  \left(  1+\sqrt{v\left(  G\right)  }\right)  $. We prove that for
all sufficiently large $n,$ there exists a graph $G=G\left(  n\right)  $ with
$\mathcal{E}\left(  G\right)  \geq n^{3/2}/2-n^{11/10}$. This implies a
conjecture of Koolen and Moulton.

We also characterize all square nonnengative matrices and all graphs with
energy close to the maximal one$.$ In particular, such graphs are quasi-random.

\textbf{Keywords: }\textit{graph energy, matrix energy, singular values,
maximal energy graphs}

\end{abstract}

Our notation is standard (e.g., see \cite{Bol98}, \cite{CDS80}, and
\cite{HoJo88}); in particular, we write $M_{m,n}$ for the set of $m\times n$
matrices, $G\left(  n\right)  $ for a graph of order $n,$ and $A^{\ast}$ for
the Hermitian adjoint of $A.$ The singular values $\sigma_{1}\left(  A\right)
\geq\sigma_{2}\left(  A\right)  \geq...$ of a matrix $A$ are the square roots
of the eigenvalues of $AA^{\ast}.$ Note that if $A\in M_{n,n}$ is a Hermitian
matrix with eigenvalues $\mu_{1}\left(  A\right)  \geq...\geq\mu_{n}\left(
A\right)  $, then $\sigma_{1}\left(  A\right)  ,...,\sigma_{n}\left(
A\right)  $ are the moduli of $\mu_{i}\left(  A\right)  $ taken in descending order.

For any $A\in M_{m,n},$ call the value $\mathcal{E}\left(  A\right)
=\sigma_{1}\left(  A\right)  +...+\sigma_{n}\left(  A\right)  $ the
\emph{energy }of $A$. Gutman \cite{Gut78}, motivated by applications in
theoretical chemistry, introduced $\mathcal{E}\left(  G\right)  =\mathcal{E}%
\left(  A\left(  G\right)  \right)  ,$ where $A\left(  G\right)  $ is the
adjacency matrix of a graph $G.$ The function $\mathcal{E}\left(  G\right)  $
has been studied intensively - see \cite{Gut05} for a survey.\ 

Recently Nikiforov \cite{Nik06} showed that if $m\leq n$ and $A\in M_{m,n}$ is
a nonnegative matrix with maximum entry $\alpha,$ then
\begin{equation}
\mathcal{E}\left(  A\right)  \leq\left(  \alpha/2\right)  \left(  m+\sqrt
{m}\right)  \sqrt{n}. \label{aub}%
\end{equation}
If in addition $\left\Vert A\right\Vert _{1}\geq n\alpha,$ then,
\begin{equation}
\mathcal{E}\left(  A\right)  \leq\left\Vert A\right\Vert _{1}/\sqrt{mn}%
+\sqrt{\left(  m-1\right)  \left(  \left\Vert A\right\Vert _{2}^{2}-\left\Vert
A\right\Vert _{1}^{2}/mn\right)  }. \label{ub}%
\end{equation}
In this note we shall investigate how tight inequality (\ref{aub}) is for
sufficiently large $m=n.$ Note that (\ref{aub}) extends an earlier bound of
Koolen and Moulton \cite{KoMo01} who proved that $\mathcal{E}\left(  G\right)
\leq\left(  n/2\right)  \left(  1+\sqrt{n}\right)  $ for any graph $G=G\left(
n\right)  ,$ and found an infinite sparse family of strongly-regular graphs
$G$ with $\mathcal{E}\left(  G\right)  =\left(  v\left(  G\right)  /2\right)
\left(  1+\sqrt{v\left(  G\right)  }\right)  $. Koolen and Moulton
\cite{KoMo01} conjectured that, for every $\varepsilon>0$, for almost all
$n\geq1,$ there exists a graph $G=G\left(  n\right)  $ with $\mathcal{E}%
\left(  G\right)  \geq\left(  1-\varepsilon\right)  \left(  n/2\right)
\left(  1+\sqrt{n}\right)  .$ We shall prove the following stronger statement.

\begin{theorem}
\label{maint}For all sufficiently large $n,$ there exists a graph $G=G\left(
n\right)  $ with $\mathcal{E}\left(  G\right)  \geq n^{3/2}/2-n^{11/10}$.
\end{theorem}

\begin{proof}
Note first that for every $G=G\left(  n\right)  ,$ we have $\sum_{i=1}%
^{n}\sigma_{i}^{2}\left(  G\right)  =tr\left(  A^{2}\left(  G\right)  \right)
=2e\left(  G\right)  $ and so
\[
2e\left(  G\right)  -\sigma_{1}^{2}\left(  G\right)  =\sigma_{2}^{2}\left(
G\right)  +...+\sigma_{n}^{2}\left(  G\right)  \leq\sigma_{2}\left(  G\right)
\left(  \mathcal{E}\left(  G\right)  -\sigma_{1}\left(  G\right)  \right)  .
\]
Hence, if $e\left(  G\right)  >0,$ then%
\begin{equation}
\mathcal{E}\left(  G\right)  \geq\sigma_{1}\left(  G\right)  +\frac{2e\left(
G\right)  -\sigma_{1}^{2}\left(  G\right)  }{\sigma_{2}\left(  G\right)
}.\label{lowb}%
\end{equation}
Let $p>11$ be a prime, $p\equiv1\left(  \operatorname{mod}\text{ }4\right)  $,
and $G_{p}$ be the Paley graph of order $p.$ Recall that $V\left(
G_{p}\right)  =\left\{  1,...,p\right\}  $ and $ij\in E\left(  G_{p}\right)  $
if and only if $i-j$ is a quadratic residue $\operatorname{mod}$ $p$. It is
known (see, e.g., \cite{Shp06}) that $G_{p}$ is a $\left(  p-1\right)
/2$-regular graph and $\sigma_{2}\left(  G_{p}\right)  =\left(  p^{1/2}%
+1\right)  /2$. Shparlinski \cite{Shp06} computed $\mathcal{E}\left(
G_{p}\right)  $ exactly, but here we need only a simple estimate. From
(\ref{lowb}) we see that
\[
\mathcal{E}\left(  G_{p}\right)  \geq\frac{p-1}{2}+\frac{p\left(  p-1\right)
/2-\left(  p-1\right)  ^{2}/4}{\left(  p^{1/2}+1\right)  /2}>\frac{p-1}%
{2}+\frac{\left(  p-1\right)  \left(  2p+1\right)  }{4\left(  p^{1/2}%
+1\right)  }>\frac{p^{3/2}}{2}.
\]
Hence, if $n$ is prime and $n\equiv1\left(  \operatorname{mod}\text{
}4\right)  $, the theorem holds. To prove it for any $n,$ recall that (see,
e.g., \cite{BHP97}, Theorem 3), for $n$ sufficiently large, there exists a
prime $p$ such that $p\equiv1\left(  \operatorname{mod}\text{ }4\right)  $ and
$p\leq n+n^{11/20+\varepsilon}.$ Suppose $n$ is large and fix some prime
$p\leq n+n^{3/5}/2.$ The average number of edges induced by a set of size $n$
in $G_{p}$ is
\[
\frac{n\left(  n-1\right)  }{p\left(  p-1\right)  }e\left(  G_{p}\right)
=\frac{n\left(  n-1\right)  }{4}.
\]
Therefore, there exists a set $X\subset V\left(  G_{p}\right)  $ with
$\left\vert X\right\vert =n$ and $e\left(  X\right)  \geq n\left(  n-1\right)
/4.$ Write $G_{n}$ for $G_{p}\left[  X\right]  $ - the graph induced by $X.$
Cauchy's interlacing theorem implies that $\sigma_{2}\left(  G_{n}\right)
\leq\sigma_{2}\left(  G_{p}\right)  $ and $\sigma_{1}\left(  G_{n}\right)
\leq\sigma_{1}\left(  G_{p}\right)  .$ Therefore, from (\ref{lowb}) we see
that%
\begin{align*}
\mathcal{E}\left(  G_{n}\right)   &  \geq\sigma_{1}\left(  G_{n}\right)
+\frac{2e\left(  G_{n}\right)  -\sigma_{1}^{2}\left(  G_{n}\right)  }%
{\sigma_{2}\left(  G_{n}\right)  }\geq\frac{\left(  n-1\right)  }{2}%
+\frac{n\left(  n-1\right)  /2-\sigma_{1}^{2}\left(  G_{p}\right)  }%
{\sigma_{2}\left(  G_{p}\right)  }\\
&  >\frac{\left(  n-1\right)  }{2}+\frac{n\left(  n-1\right)  /2-\left(
n+n^{3/5}/2\right)  ^{2}/4}{\left(  \sqrt{n+n^{3/5}/2}+1\right)  /2}%
>\frac{n^{3/2}}{2}-n^{11/10},
\end{align*}
completing the proof.
\end{proof}

Clearly Theorem \ref{maint} implies that inequality (\ref{aub}) is tight for
all $n\times n$ nonnegative matrices as well. To the end of the note we shall
characterize all square nonnegative matrices and all graphs with energy close
to the maximal one.

\begin{theorem}
For every $\varepsilon>0$ there exists $\delta>0$ such that, for all
sufficiently large $n$, if $A\in M_{n,n}$ is a nonnegative matrix with maximum
entry $\alpha>2\varepsilon,$ and $\mathcal{E}\left(  A\right)  \geq
\alpha\left(  1/2-\delta\right)  n^{3/2},$ \ then the following conditions hold:

\qquad(i) $a_{ij}>\left(  1-\varepsilon\right)  \alpha$ for at least $\left(
1/2-\varepsilon\right)  n^{2}$ entries $a_{ij}$ of $A;$

\qquad(ii) $a_{ij}<\varepsilon\alpha$ for at least $\left(  1/2-\varepsilon
\right)  n^{2}$ entries $a_{ij}$ of $A;$

\qquad(iii) $\left\vert \sigma_{1}\left(  A\right)  -\alpha n/2\right\vert
<\varepsilon\alpha n;$

\qquad(iv) $\sigma_{2}\left(  A\right)  <\varepsilon\alpha n;$

\qquad(v) $\left\vert \sigma_{i}\left(  A\right)  -\alpha n^{1/2}/2\right\vert
<\varepsilon\alpha n^{1/2}$ for all $\varepsilon n\leq i\leq\left(
1-\varepsilon\right)  n.$
\end{theorem}

\begin{proof}
Without loss of generality we shall assume that $\alpha=1.$ Note that
$\mathcal{E}\left(  A\right)  \leq\sqrt{n\left\Vert A\right\Vert _{2}^{2}%
}<\sqrt{n\left\Vert A\right\Vert _{1}}$ and so $\left\Vert A\right\Vert
_{1}>n.$ Summarizing the essential steps in the proof of (\ref{aub}) (see
\cite{Nik06}), we have
\begin{align*}
\left(  1/2-\delta\right)  n^{3/2}  &  \leq\mathcal{E}\left(  A\right)
\leq\sigma_{1}\left(  A\right)  +\sqrt{\left(  n-1\right)  \left(  \left\Vert
A\right\Vert _{2}^{2}-\sigma_{1}^{2}\left(  A\right)  \right)  }\\
&  \leq\left\Vert A\right\Vert _{1}/n+\sqrt{\left(  n-1\right)  \left(
\left\Vert A\right\Vert _{2}^{2}-\left\Vert A\right\Vert _{1}^{2}%
/n^{2}\right)  }\\
&  \leq\left\Vert A\right\Vert _{1}/n+\sqrt{\left(  n-1\right)  \left(
\left\Vert A\right\Vert _{1}-\left\Vert A\right\Vert _{1}^{2}/n^{2}\right)
}\\
&  \leq\frac{n}{2}\left(  \sqrt{n}+1\right)  .
\end{align*}
From continuity, we see that \emph{(i), (ii), (iii),} and \emph{(v)} hold for
$\delta$ small and $n$ large enough. Noting that
\[
\left(  1/2-\delta\right)  n^{3/2}\leq\mathcal{E}\left(  A\right)  \leq
\sigma_{1}\left(  A\right)  +\sigma_{2}\left(  A\right)  +\sqrt{\left(
n-2\right)  \left(  \left\Vert A\right\Vert _{1}-\sigma_{1}^{2}\left(
A\right)  -\sigma_{1}^{2}\left(  A\right)  \right)  },
\]
in view of \emph{(i), (ii), }and\emph{ (iii)},\emph{ }we see that \emph{(iv)}
holds as well for $\delta$ small and $n$ large enough, completing the proof.
\end{proof}

Note that the characterization given by \emph{(i) - (v) }is complete, because
if \emph{(v)} alone holds, then
\[
\mathcal{E}\left(  A\right)  >\alpha\left(  1-2\varepsilon\right)  n\left(
\frac{1}{2}-\varepsilon\right)  n^{1/2}>\alpha\left(  1-2\varepsilon\right)
n^{3/2}.
\]
Also \emph{(i) - (iv)} imply that any graph whose energy is close to
$n^{3/2}/2$ is quasi-random in the sense of \cite{CGW89}.

We conclude with the following interesting statement, whose proof is omitted.

\begin{theorem}
For every $\varepsilon>0$ there exists $\delta>0$ such that, if $n$ is large
enough$,$ $A\in M_{n,n}$ is a nonnegative matrix with maximum entry
$0<\varepsilon<1/2,$ and $\mathcal{E}\left(  A\right)  \geq\left(
1/2-\delta\right)  n^{3/2},$ then $\mathcal{E}\left(  E_{n}-A\right)
\geq\left(  1/2-\varepsilon\right)  n^{3/2},$ where $E_{n}\in M_{n,n}$ is the
matrix of all ones.
\end{theorem}

It would be interesting to determine how tight is inequality (\ref{aub}) for
nonsquare matrices.


\begin{thebibliography}{99}                                                                                               %


\bibitem {BHP97}R. C. Baker, G. Harman, J. Pintz, The exceptional set for
Goldbach's problem in short intervals, \emph{Sieve methods, exponential sums,
and their applications in number theory }(Cardiff, 1995),\ pp. 1--54, LMS
Lecture Notes Ser., 237, Cambridge Univ. Press, Cambridge, 1997.

\bibitem {Bol98}B. Bollob\'{a}s, \emph{Modern Graph Theory}\textit{,} Graduate
Texts in Mathematics, 184, Springer-Verlag, New York (1998), xiv+394 pp.

\bibitem {CDS80}D. Cvetkovi\'{c}, M. Doob, and H. Sachs, \emph{Spectra of
Graphs,} VEB Deutscher Verlag der Wissenschaften, Berlin, 1980, 368 pp.

\bibitem {CGW89}F. R. K. Chung, R. L. Graham, and R. M. Wilson, Quasi-random
graphs, Combinatorica 9(1989), 345--362.

\bibitem {Gut78}I. Gutman, The energy of a graph, Ber. Math.-Stat. Sekt.
Forschungszent. Graz 103 (1978), 1--22.

\bibitem {Gut05}I. Gutman, Topology and stability of conjugated hydrocarbons.
The dependence of total pi-electron energy on molecular topology, J. Serb.
Chem. Soc. 70 (2005) 441--456.

\bibitem {HoJo88}R. Horn and C. Johnson, \emph{Matrix Analysis,} Cambridge
University Press, Cambridge, 1985, xiii+561 pp.

\bibitem {KoMo01}J.H. Koolen, V. Moulton, Maximal energy graphs, Adv. Appl.
Math. 26 (2001), 47--52.

\bibitem {Nik06}V. Nikiforov, The energy of graphs and matrices, to appear in
J. Math. Anal. Appl.

\bibitem {Shp06}I. Shparlinski, On the energy of some circulant graphs, Linear
Algebra Appl. 414 (2006), 378--382.
\end{thebibliography}
\end{document}